# Upper and Lower Bounds on the Number of Faults a System Can Withstand Without Repairs*

Michel Goemans†    Nancy Lynch‡    Isaac Saias§

Laboratory for Computer Science
Massachusetts Institute of Technology
Cambridge, MA  02139

ABSTRACT
We consider the following scheduling problem. A system is composed of $n$ processors drawn from a pool of $N$. The processors can become faulty while in operation and faulty processors never recover. A report is issued whenever a fault occurs. This report states only the existence of a fault, but does not indicate its location. Based on this report, the scheduler can reconfigure the system and choose another set of $n$ processors. The system operates satisfactorily as long as at most $f$ of the $n$ selected processors are faulty. We exhibit a scheduling strategy allowing the system to operate satisfactorily until approximately $(N/n)f$ faults are reported in the worst case. Our precise bound is tight.

## 1  Introduction

Many control systems are subject to failures that can have dramatic effects. One simple way to deal with this problem is to build in some redundancy so that the whole system is able to function even if parts of it fail. In a general situation, the system's manager has access to some observations allowing it to control the system efficiently. Such observations bring information about the state of the system that might consist of partial fault reports. The available controls might include repairs and/or replacement of faulty processors.

To model the problem, one needs to make assumptions regarding the occurrence of faults. Typically, they are assumed to occur according to some stochastic process. To make the model more tractable, one often considers the process to be memoryless, i.e. faults occur according to some exponential distribution. However, to be more realistic, many complica-tions and variations can be introduced in the stochastic model, and they complicate the time analysis. Examples are: a processor might become faulty at any time or only during specific operations; the fault rate might vary according to the work load; faults might occur independently among the processors or may depend on proximity. The variations seem endless and the results are rarely general enough so as to carry some information or methodology from one model to another.

---

*Research supported by research contracts AFOSR-89-0271, ONR-N00014-91-J-1046, NSF-CCR-8915206 and DARPA-N00014-89-J-1988.
†e-mail: goemans@math.mit.edu
‡e-mail: lynch@theory.mit.edu
§e-mail: saias@theory.lcs.mit.edu

One way to derive general results, independent of the specific assumptions about the time of occurrence of faults, is to adopt a *logical time*, that, instead of following an absolute frame, is incremented only at each occurrence of a fault. Within this framework, we measure the maximal number of faults to be observed until the occurrence of a crash instead of the maximal time of survival of a system until the occurrence of a crash.

As an introduction to this general situation, we make the following assumptions and simplifications:

**Redundancy of the system:** We assume the existence of a pool $\mathcal{N}$ composed of $N$ identical processors from among which, at every time $t$, a set $S_t$ of $n$ processors is selected to configure the system. The system works satisfactorily as long as at least $n - f$ processors among the $n$ currently in operation are not faulty. tolerate more than $f$ faults at any given time: it stops functioning if $f + 1$ processors among these $n$ processors are faulty.

**Occurrence of faults, reports and logical time:** We consider the situation in which failures do not occur simultaneously and where, whenever a processor fails, a report is issued, stating that a failure has occurred, but without specifying the location of the failure. (Reporting additional information might be too expensive or time consuming.) Based on these reports, the scheduler might decide to reconfigure the system whenever such failure is reported. As a result, we restrict our attention to the discrete model, in which time $t$ corresponds to the $t$-th failure in the system.

**Repairs:** No repair is being performed.

**Deterministic Algorithms:** We assume that the scheduler does not use randomness.

Since the universe consists of only $N$ processors, and one processor fails at each time, no scheduling policy can guarantee that the system survives beyond time $N$. (A better a priori upper bound is $N - n + f + 1$: at this time, only $n - f - 1$ processors are still non-faulty. This does not allow for the required quorum of $n - f$ non-faulty processors.) But some scheduling policies seem to allow the system to survive longer than others. An obviously bad policy is to choose $n$ processors once and for all and never to change them: the system would then collapse at time $f + 1$. This paper investigates the problem of determining the best survival time.

This best survival time is defined from a worst-case point-of-view: a *given* scheduler allows the system to survive (up to a certain time) only if it allows it to survive against *all* possible failure patterns in which one processor fails at each time.

Our informal description so far apparently constrains the faults to occur in on-line fashion: for each $t$, the $t$-th fault occurs before the scheduler decides the set $S_{t+1}$ to be used subsequently. However, since we have assumed that no reports about the locations of the faults are available, there is no loss of generality in requiring the sets $S_t$ to be determined a priori. (Of course, in practice, some more precise fault information may be available, and each set $S_t$ would depend on the fault pattern up to time $t$.) Also, as we have assumed a *deterministic* scheduler, we can assume that the decisions $S_1, \ldots, S_N$ are *revealed* before the occurrence of any fault. We express this by saying that the faults occur in an off-line fashion.

## 2  The Model

Throughout this paper, we fix a universal set $\mathcal{N}$ of processors, and let $N$ denote its cardinality. We also fix a positive integer $n$ ($n \leq N$) representing the number of processors that are needed at each time period, and a positive integer $f$ representing the number of failures that can be tolerated ($f < n$).

We model the situation described in the introduction as a simple game between two entities, a *scheduler* and an *adversary*. The game consists of only one round, in which the scheduler plays first and the adversary second. The scheduler plays by selecting a sequence of $N$ sets of processors (the *schedule*), each set of size $n$, and the adversary responds by choosing, from each set selected by the scheduler, a processor to kill. We consider only sequences of size $N$ because the system must collapse by time $N$, since, at each time period, a new processor breaks down.

Formally, a schedule $\mathcal{S}$ is defined to be a finite sequence, $S_1, \ldots, S_N$, of subsets of $\mathcal{N}$, such that $|S_t| = n$ for all $t$, $1 \leq t \leq N$. An $\mathcal{S}$-*adversary* $\mathcal{A}$ is defined to be a finite sequence, $s_1, \ldots, s_N$, of elements of $\mathcal{N}$ such that $s_t \in S_t$ for every $t$.

Now let $\mathcal{S}$ be a schedule, and $\mathcal{A}$ an $\mathcal{S}$-adversary. Define the *survival time*, $T(\mathcal{S}, \mathcal{A})$, to be the largest value of $t$ such that, for all $u \leq t$, $|\{s_1, \ldots s_u\} \cap S_u| \leq f$. That is, for all time periods $u$ up to and including time period $t$, there are no more than $f$ processors in the set $S_u$ that have failed by time $u$.

We are interested in the minimum survival time for a particular schedule, with respect to arbitrary adversaries. Thus, we define the *minimum survival time* for a schedule, $T(\mathcal{S})$, to be $T(\mathcal{S}) \stackrel{\text{def}}{=} \min_{\mathcal{A}} T(\mathcal{S}, \mathcal{A})$. In this definition, the minimum is taken over all $\mathcal{S}$-adversaries. An adversary $\mathcal{A}$ for which $T(\mathcal{S}) = \min_{\mathcal{A}} T(\mathcal{S}, \mathcal{A})$ is said to be *minimal* for $\mathcal{S}$. Finally, we are interested in determining the schedule that guarantees the greatest minimum survival time. Thus, we define the *optimum survival time* $T_{\text{opt}}$, to be $\max_{\mathcal{S}} T(\mathcal{S}) = \max_{\mathcal{S}} \min_{\mathcal{A}} T(\mathcal{S}, \mathcal{A})$. Also define a schedule $\mathcal{S}$ to be *optimum* provided that $T(\mathcal{S}) = T_{\text{opt}}$. Our objectives in this paper are to compute $T_{\text{opt}}$ as a function of $N$, $n$ and $f$, to exhibit an optimum schedule, and to determine a minimal adversary for each schedule.

## 3  The Result

Recall that $1 \leq f < n \leq N$ are three fixed integers. Our main result is stated in terms of the following function defined on the set of positive real numbers (see Figure 1):

$$h_{n,f}(k) \stackrel{\text{def}}{=} \left\lfloor \frac{k}{n} \right\rfloor f + \left( k - \left\lfloor \frac{k}{n} \right\rfloor n + f - n \right)^+,$$

where $(x)^+ = \max(x, 0)$. In particular, $h_{n,f}(k) = \frac{k}{n}f$ when $n$ divides $k$.

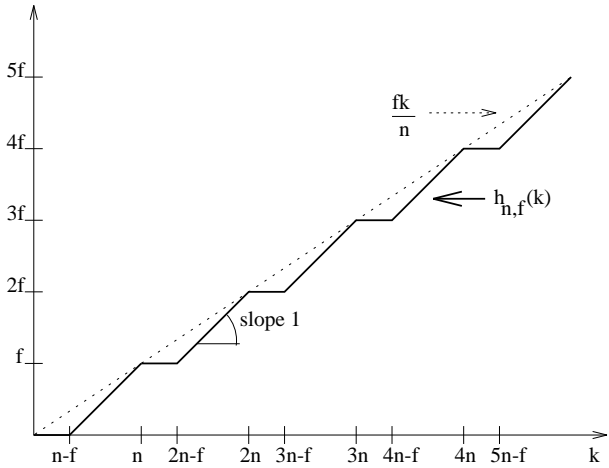

Figure 1. The function $h_{n,f}(k)$

The main result of this paper is:

**Theorem 3.1**
$$T_{\text{opt}} = h_{n,f}(N).$$

We will present our proof in two lemmas proving respectively that $T_{\text{opt}}$ is no smaller and no bigger than $h_{n,f}(N)$.

**Lemma 3.2**
$$T_{\text{opt}} \geq h_{n,f}(N).$$

*Proof:* Consider the schedule $\mathcal{S}_{\text{trivial}}$ in which the $N$ processors are partitioned into $\lfloor \frac{N}{n} \rfloor$ batches of $n$ processors each and one batch of $p = N - \lfloor \frac{N}{n} \rfloor n$. Each of the first $\lfloor \frac{N}{n} \rfloor$ batches is used $f$ time periods and then set aside. Then, the last batch of processors along with any $n - p$ of the processors set aside is used for $(f + p - n)^+$ time periods. It is easy to see that no adversary can succeed in killing $f + 1$ processors within a batch before this schedule expires. ∎

In order to prove the other direction of Theorem 3.1, we need the following result about the rate of increase of the function $h_{n,f}(k)$.

**Lemma 3.3** *For $0 \leq k$ and $0 \leq l \leq n$ we have $h_{n,f}(k) \leq h_{n,f}(k+l) + n - l - f$.*

*Proof:* Notice first that $h_{n,f}(k) = h_{n,f}(k+n) - f$ for all $k \geq 0$. Moreover, the function $h$ increases at a sublinear rate (see Figure 1) so that, for $p, q \geq 0$, we have $h_{n,f}(p+q) \leq h_{n,f}(p) + q$. Letting $p = k + l$ and $q = n - l$, we obtain

$$h_{n,f}(k) = h_{n,f}(k+n) - f \leq h_{n,f}(k+l) + n - l - f,$$

which proves the lemma.    ∎

## 4   The Upper Bound

In this section we establish the other direction of the main theorem. We begin with some general graph theoretical definitions.

**Definition 4.1**

- For every vertex $v$ of a graph $G$, we let $\gamma_G(v)$ denote the set of vertices adjacent to $v$. We can extend this notation to sets: for all sets $C$ of vertices $\gamma_G(C) \stackrel{\text{def}}{=} \cup_{v \in C} \gamma_G(v)$.

- For every bipartite graph $G$, $\nu(G)$ denotes the size of a maximum matching of $G$.

- For every pair of positive integers $L, R$, a *left totally ordered* bipartite graph of size $(L, R)$ is a bipartite graph with bipartition $\mathcal{L}, \mathcal{R}$, where $\mathcal{L}$ is a totally ordered set of size $L$ and $\mathcal{R}$ is a set of size $R$. We label $\mathcal{L} = \{a_1, \ldots, a_L\}$ so that, $a_i < a_j$ for every $1 \leq i < j \leq L$. For every $\mathcal{L}' \subseteq \mathcal{L}$ and $\mathcal{R}' \subseteq \mathcal{R}$, the subgraph induced by $\mathcal{L}'$ and $\mathcal{R}'$ is a left totally ordered bipartite graph with the total order on $\mathcal{L}$ inducing the total order on $\mathcal{L}'$.

- Let $G$ be a left totally ordered bipartite graph of size $(L, R)$. For $t = 1, \ldots, L$, we let $I_t(G)$ denote the left totally ordered subgraph of $G$ induced by the subsets $\{a_1, a_2, \ldots, a_{t-1}\} \subseteq \mathcal{L}$ and $\gamma_G(a_t) \subseteq \mathcal{R}$.

Let us justify quickly the notion of left total order. In this definition, we have in mind that $\mathcal{L}$ represents the labels attached to the different times, and that $\mathcal{R}$ represents the labels attached to the available processors. The times are naturally ordered. The main argument used in the proof is to reduce an existing schedule to a shorter one. In doing so, we in particular select a subsequence of times. Although these times are not necessarily consecutive, they are still naturally ordered. The total order on $\mathcal{L}$ is the precise notion formalizing the ordering structure characterizing time.

Consider a finite schedule $S = S_1, \ldots, S_T$. In graph theoretic terms, it can be represented as a left totally ordered bipartite graph $G$ with bipartition $\mathcal{T} = \{1, 2, \ldots, T\}$ and $\mathcal{N} = \{1, 2, \ldots, N\}$. There is an edge between vertex $t \in \mathcal{T}$ and vertex $i \in \mathcal{N}$ if the processor $i$ is selected at time $t$. The fact that, for all $t$, $|S_t| = n$ translates into the fact that vertex $t \in \mathcal{T}$ has degree $n$. For such a bipartite graph, the game of the adversary consists in selecting one edge incident to each vertex $t \in \mathcal{T}$.

Observe that the adversary can kill the schedule at time $t$ if it has already killed, before time $t$, $f$ of the $n$ processors used at time $t$. It then kills another one at time $t$ and the system collapses. In terms of the graph $G$, there exists an adversary that kills the schedule at time $t$ if and only if the subgraph $I_t(G)$ has a matching of size $f$, i.e. $\nu(I_t(G)) \geq f$. Therefore, the set $\mathcal{P}$ that we now define represents the set of integers $L$ and $R$ for which there exists a schedule that survives at time $L$, when $R$ processors are available.

**Definition 4.2** Let $L$ and $R$ be two positive integers. $(L, R) \in \mathcal{P}$ iff there exists a left totally ordered bipartite graph $G$ of size $(L, R)$ with bipartition $\mathcal{L}$ and $\mathcal{R}$ satisfying the two following properties:

1. All vertices in $\mathcal{L}$ have degree exactly equal to $n$,

2. For every $t = 1, \ldots, |\mathcal{L}|$, all matchings in $I_t(G)$ have size at most equal to $f - 1$, i.e. $\nu(I_t(G)) \leq f - 1$.

The main tool used in the proof of Theorem 3.1 is the following duality result for the maximum bipartite matching problem, known as Ore's Deficiency Theorem [3]. A simple proof of this theorem and related results can be found in [2].

**Theorem 4.1** *Let $G$ be a bipartite graph with bipartition $A$ and $B$. Then the size $\nu(G)$ of a maximum matching is given by the formula:*

$$\nu(G) = \min_{C \subseteq B} \left[ |B - C| + |\gamma_G(C)| \right]. \tag{1}$$

The following lemma is crucial for our proof.

**Lemma 4.2** *There are no positive integers $L$ and $R$ such that $(L, R) \in \mathcal{P}$ and such that $L > h_{n,f}(R)$.*

*Proof:*

Working by contradiction, consider two positive integers $L$ and $R$ such that $(L, R) \in \mathcal{P}$ and $L > h_{n,f}(R)$. We first show the existence of two integers $L'$ and $R'$ such that $L' < L$, $(L', R') \in \mathcal{P}$ and $L' > h_{n,f}(R')$.

Let $\mathcal{L} = \{a_1, a_2, \ldots, a_L\}$ and $\mathcal{R} = \{b_1, b_2, \ldots, b_R\}$ be the bipartition of the graph $G$ whose existence is ensured by the hypothesis $(L, R) \in \mathcal{P}$.

We apply Theorem 4.1 to the graph $I_L(G)$ where we set $A = \{a_1, a_2, \ldots, a_{L-1}\}$ and $B = \gamma_G(a_L)$. Let $C$ denote a subset of $B$ for which the minimum in (1) is attained. ($C$ is possibly empty.) Define $\mathcal{L}' \stackrel{\text{def}}{=} \mathcal{L} - (\{a_L\} \cup \gamma_{I_L(G)}(C))$ and $\mathcal{R}' \stackrel{\text{def}}{=} \mathcal{R} - C$ and let $L'$ and $R'$ denote the cardinalities of $\mathcal{L}'$ and $\mathcal{R}'$. Hence, $L' = L - 1 - |\gamma_{I_L(G)}(C)|$ so that $L' < L$. Consider the bipartite subgraph $G'$ of $G$ induced by the set of vertices $\mathcal{L}' \cup \mathcal{R}'$. In other words, in order to construct $G'$ from $G$, we remove the set $C \cup \{a_L\}$ of vertices and *all* vertices adjacent to some vertex in $C$. We have illustrated this construction in Figure 2. In that specific example, $n = 4$, $f = 3$, $L = 6$ and $R = 7$, while $h_{4,3}(7) = 5$. One can show that $C = \{b_5, b_6, b_7\}$ and as a result $G'$ is the graph induced by the vertices $\{a_1, a_2, a_3, a_4, b_1, b_2, b_3, b_4\}$. The graph $G'$ has size $(L', R') = (4, 4)$.

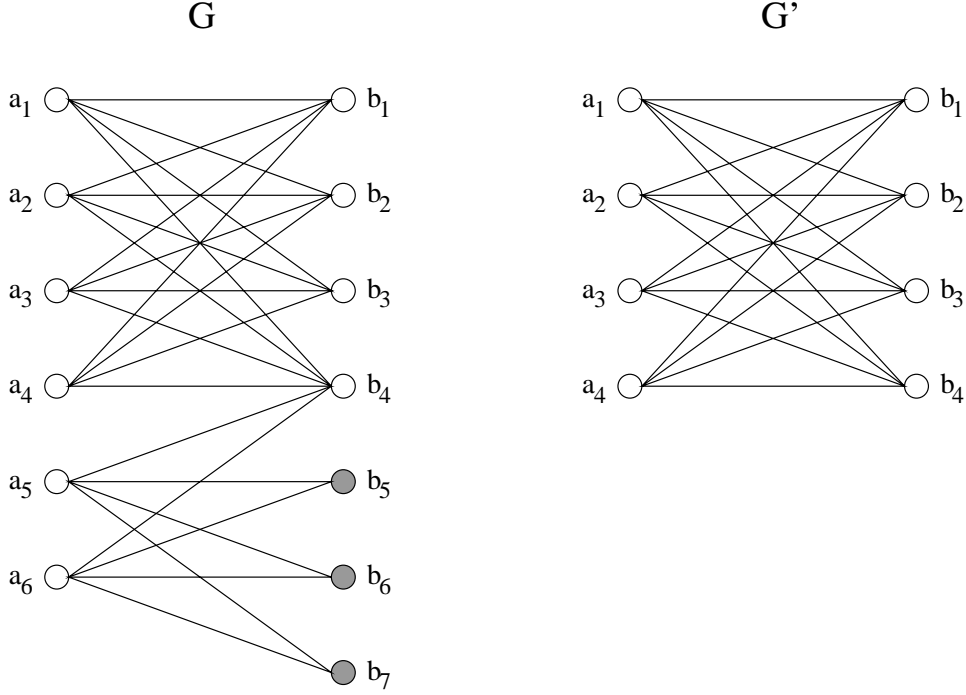

Figure 2. An example of the construction of $G'$ from $G$. The vertices in $C$ are darkened.

We first show that $(L', R') \in \mathcal{P}$. Since the vertices in $\mathcal{L}'$ correspond to the vertices of $\mathcal{L} - \{a_L\}$ not connected to $C$, their degree in $G'$ is also $n$. Furthermore, $G'$, being a subgraph of $G$, inherits property 2 of Definition 4.2. Indeed, assume that there is a vertex $a_{t'}$ in $G'$ such that $I_{t'}(G')$ has a matching of size $f$. Let $t$ be the label of the corresponding vertex in graph $G$. Since the total order on $\mathcal{L}'$ is induced by the total order on $\mathcal{L}$, $I_{t'}(G')$ is a subgraph of $I_t(G)$. Therefore, $I_t(G)$ would also have a matching of size $f$, a contradiction.

Let us show that $L' > h_{n,f}(R')$. The assumption $(L, R) \in \mathcal{P}$ implies that $f - 1 \geq \nu(I_L(G))$. Using Theorem 4.1 and the fact that $B = \gamma_G(L)$ has cardinality $n$, this can be rewritten as

$$\begin{aligned} f - 1 &\geq \nu(I_L(G)) = |B - C| + |\gamma_{I_L(G)}(C)| \\ &= n - |C| + |\gamma_{I_L(G)}(C)|. \end{aligned} \quad (2)$$

Since $C \subseteq B \subseteq \mathcal{R}$, we have that $0 \leq |C| \leq n \leq R$ and, thus, the hypotheses of Lemma 3.3 are satisfied for $k = R - |C|$ and $l = |C|$. Therefore, we derive from the lemma that

$$h_{n,f}(R') = h_{n,f}(R - |C|) \leq h_{n,f}(R) + n - |C| - f.$$

Using (2), this implies that

$$h_{n,f}(R') \leq h_{n,f}(R) - |\gamma_{I_L(G)}(C)| - 1.$$

By assumption, $L$ is strictly greater than $h_{n,f}(R)$, implying

$$h_{n,f}(R') < L - 1 - |\gamma_{I_L(G)}(C)|.$$

But the right-hand-side of this inequality is precisely $L'$, implying that $L' > h_{n,f}(R')$.

We have therefore established that for all integers $L$ and $R$ such that $(L, R) \in \mathcal{P}$ and $L > h_{n,f}(R)$, there exists two integers $L'$ and $R'$ such that $L' < L$, $(L', R') \in \mathcal{P}$ and $L' > h_{n,f}(R')$. Among all such pairs $(L, R)$, we select the pair for which $L$ is minimum. By the result that we just established, we obtain a pair $(L', R')$ such that $(L', R') \in \mathcal{P}$ and $L' < L$. This contradicts the minimality of $L$.

∎

**Lemma 4.3**
$$T_{\text{opt}} \leq h_{n,f}(N).$$

*Proof:* By assumption, $(T_{\text{opt}}, N) \in \mathcal{P}$. Hence this result is a direct consequence of Lemma 4.2 .

∎

This Lemma along with Lemma 3.2 proves Theorem 3.1.

In the process of proving Lemma 3.2 we proved that $S_{\text{trivial}}$ is an optimum schedule. On the other hand, the interpretation of the problem as a graph problem also demonstrates that the adversary has a polynomial time algorithm for finding an optimum killing sequence for each schedule $S$. When provided with $S$, the adversary needs only to compute a polynomial number (actually fewer than $N$) of maximum bipartite matchings, for which well known polynomial algorithms exist (for the fastest known, see [1]).

## 5  Future Research

The problem solved in this paper is a first step towards modeling complex resilient systems and there are many interesting extensions. We mention only a few.

An interesting extension is to consider the case of a system built up of processors of different types. For instance consider the case of a system built up of a total of $n$ processors, that is reconfigured at each time period and that needs at least $g_1$ non-faulty processors of type 1 and at least $g_2$ non-faulty processors of type 2 in order to function satisfactorily. Assume also that these processors are drawn from a pool $\mathcal{N}_1$ of $N_1$ processors of type 1 and a pool $\mathcal{N}_2$ of $N_2$ processors of type 2, that $\mathcal{N}_1 \cap \mathcal{N}_2 = \emptyset$, that that there are no repairs. It is easy to see that the optimum survival time $T_{\text{opt}}$ is at least the survival time of every strategy for which the number of processors of type 1 and type 2 is kept constant throughout. Hence:

$$T_{\text{opt}} \geq \max_{\{(n_1, n_2); n_1 + n_2 = n\}} \min\left(h_{n_1, n_1 - g_1}(N_1), h_{n_2, n_2 - g_2}(N_2)\right).$$

It would be an interesting question whether $T_{\text{opt}}$ is exactly equal to this value or very close to it.

Extend the definition of a *scheduler* to represent a randomized scheduling protocol. (Phrased in this context, the result presented in this paper is only about deterministic scheduling protocols.) A scheduler is called *adversary-oblivious* if it decides the schedule independently of the choices $s_1, s_2, \ldots$ made by the adversary. An *off-line* adversary is an adversary that has access to the knowledge of the full schedule $S_1, S_2, \ldots$ before deciding the full sequence $s_1, s_2, \ldots$ Note that, by definition, off-line adversaries make sense only with adversary-oblivious schedulers. By comparison, an *on-line* adversary decides for each time $t$ which processor $s_t$ to kill, without knowing the future schedule: at each time $t$ the adversary decides $s_t$ based on the sole knowledge of $S_1, \ldots, S_t$ and of $s_1, \ldots, s_{t-1}$. In this more general framework, the quantity we want to determine is

$$T_{\text{opt}} \stackrel{\text{def}}{=} \max_{\mathcal{S}} \min_{\mathcal{A}} E\left[T(\mathcal{S}, \mathcal{A})\right]. \tag{3}$$

For an adversary-oblivious, randomized scheduler, one can consider two cases based on whether the adversary is on-line or off-line. As is easily seen, if the adversary is off-line, randomness does not help in the design of optimal schedulers: introducing randomness in the schedules cannot increase the survival time if the adversary gets full knowledge of the schedule before committing to any of its choices. As a result, the off-line version corresponds to the situation investigated in this paper.

It would be of interest to study the online version of Problem (3). On-line adversaries model somewhat more accurately practical situations: faults naturally occur in an on-line fashion and the role of the program designer is then to design a scheduler whose expected performance is optimum. Hence, comparing the two versions of Problem 3 would allow to understand how much randomness can help in the design of optimum, adversary-oblivious, schedulers.

For instance, in the case where $N = 4, n = 2$ and $f = 1$, and where on-line adversaries are considered, a case analysis shows that $T_{\text{opt}}$ is equal to $9/4$ for randomized algorithms. A direct application of Theorem 3.1 shows that that $T_{\text{opt}} = 2$ for deterministic algorithms.

Going towards even more realistic and complex situations, we can also study the case where the scheduler is provided at each time with some partial information about the fault process.

*Acknowledgements:* We want to thank Stuart Adams from the Draper Laboratories who suggested the problem to us and George Varghese for helpful suggestions.

## 6   REFERENCES

[1] J.E. Hopcroft and R.M. Karp, "A $n^{5/2}$ algorithm for maximum matching in bipartite graphs", *SIAM J. Computing*, 2, 225–231 (1973).

[2] L. Lovász and M. Plummer, "Matching theory", *North-Holland Mathematical Studies 121*, Annals of Discrete Mathematics, 29, (1986).

[3] O. Ore, "Graphs and matching theorems", *Duke Math. J.*, 22, 625–639 (1955).